\documentclass[a4paper]{amsart}

\usepackage{amssymb}

\newenvironment{pf}{\begin{proof}}{\end{proof}}

\renewcommand{\mod}{\operatorname{mod}}
\DeclareMathOperator{\CM}{CM}

\DeclareMathOperator{\Gal}{Gal}
\newcommand{\spitz}[1]{\langle #1\rangle}

\DeclareMathOperator{\Hom}{Hom}
\DeclareMathOperator{\End}{End}

\DeclareMathOperator{\rad}{rad}

\DeclareMathOperator{\Aut}{Aut}
\DeclareMathOperator{\Inn}{Inn}

\DeclareMathOperator{\ch}{char}

\newcommand{\quat}[3]{\bigl(\frac{#1,\, #2}{#3}\bigr)}
\newcommand{\quatt}[3]{\bigl(\frac{#1,\, #2}{#3}\bigr]}

\DeclareMathOperator{\matring}{M}

\numberwithin{equation}{section}

\theoremstyle{plain}
\newtheorem{Prop}{Proposition}[section]
\newtheorem{Thm}[Prop]{Theorem}
\newtheorem{Lem}[Prop]{Lemma}
\newtheorem{Cor}[Prop]{Corollary}

\theoremstyle{definition}
\newtheorem{Defi}[Prop]{Definition}
\newtheorem{Ex}[Prop]{Example}
\newtheorem{Exs}[Prop]{Examples}
\newtheorem{Rem}[Prop]{Remark}
\newtheorem{Rems}[Prop]{Remarks}
\newtheorem*{Ack}{Acknowledgements}

\begin{document}

\title[Parameter curves]{Parameter curves for the regular
  representations of tame bimodules}


\author{Dirk Kussin}

\address{Institut f\"ur Mathematik\\
Universit\"at Paderborn\\
33095 Paderborn\\
Germany}

\email{dirk@math.upb.de}


\begin{abstract}
  We present results and examples which show that the consideration of
  a certain tubular mutation is advantageous in the study of
  noncommutative curves which parametrize the simple regular
  representations of a tame bimodule. We classify all tame bimodules
  where such a curve is actually commutative, or in different words,
  where the unique generic module has a commutative endomorphism ring.
  This extends results from~\cite{kussin:07} to arbitrary
  characteristic; in characteristic two additionally inseparable cases
  occur. Further results are concerned with autoequivalences fixing
  all objects but not isomorphic to the identity functor.
\end{abstract}

\maketitle

\section{Introduction}

The notion of tameness for finite dimensional algebras defined over a
field $k$ which is not algebraically closed is not well understood.  A
study of the class of tame hereditary algebras is indispensable for
understanding tameness in general, since these algebras are the
easiest tame algebras in the sense that there is precisely one generic
module~\cite{ringel:79,crawleyboevey:91b}.

Let $\Lambda$ be a tame hereditary algebra over a field $k$. A
commutative curve which parametrizes the simple regular
$\Lambda$-modules was studied in~\cite{crawleyboevey:91}.
In~\cite{bgl:87} the authors introduced a noncommutative curve
$\mathbb{X}$, also parametrizing the simple regular $\Lambda$-modules
and whose centre is the curve from~\cite{crawleyboevey:91}. This
noncommutative curve contains the full information on $\Lambda$, its
function (skew) field $k(\mathbb{X})$ is the endomorphism ring of the
unique generic module, and it is defined by the so-called
preprojective algebra, which is an orbit algebra (see
Section~\ref{sec:efficient-auto}) formed with respect to the inverse
Auslander-Reiten translation $\tau^{-}$. For the study of this curve
$\mathbb{X}$ it is advantageous to replace $\tau^{-}$ (whenever
possible) by a certain tubular mutation and to form an orbit algebra
with respect to such a functor. We refer to~\cite{kussin:07} for the
foundations. 

By the technique of reducing and inserting
weights~\cite{geiglelenzing:91,lenzing:97,kussin:07} it is essentially
sufficient to study the underlying tame bimodule; in this case all
tubes are homogeneous. Thus, in the terminology of~\cite{kussin:07} we
study the noncommutative homogeneous curves of genus zero, or in the
terminology of~\cite{lenzing:97,kussin:07} the homogeneous exceptional
curves. We exploit our orbit algebra in particular
\begin{itemize}
\item to classify all tame bimodules (arbitrary characteristic) and
  the corresponding function fields which yield a \emph{commutative}
  curve $\mathbb{X}$; in particular we exhibit an inseparable example
  where $\mathbb{X}$ is a commutative curve but not a Brauer-Severi
  curve;
\item to get precise information on the so-called ghost group
  $\mathcal{G}$, and accordingly also on the Auslander-Reiten
  translation as functor;
\item to describe $\mathbb{X}$ and its function field $k(\mathbb{X})$
  explicitly in many examples, in particular for all tame bimodules
  over a finite field (characteristic different from $2$).
\end{itemize}
The results on the commutative cases include in particular the
well-known classification of the (commutative) function fields in one
variable of genus zero in the classical sense (we refer
to~\cite[Ch.~16.]{artin:67},
\cite[Thm.~1.1]{buchweitz_eisenbud_herzog:87}, \cite{lang_tate:52}).
Our noncommutative treatment sheds new light on this classical topic.
We stress that in our setting also function fields are allowed which
are not separably generated over $k$.

\section{Efficient automorphisms}\label{sec:efficient-auto}

For the study of a parametrization of the indecomposable modules over
a tame hereditary algebra it is essentially sufficient to consider the
bimodule case, which is also called the homogeneous case. Thus, let
$M={}_{F}M_{G}$ be a tame bimodule ($F$ and $G$ division algebras) of
finite dimension over the centrally acting field $k$, and let
$\Lambda={\footnotesize
\begin{pmatrix}
  G & 0\\
  M & F
\end{pmatrix}}$ be the corresponding tame hereditary (bimodule)
algebra. Recall that \emph{tame} means $\dim {}_F M \cdot\dim M_G =4$.
We refer to~\cite{dlabringel:76,ringel:76} for further information on
bimodules and their representation theory. We will assume throughout
that (without loss of generality) $k$ is the centre of $M$.

The geometric structure of $\mathbb{X}$, the index set of the simple
regular $\Lambda$-modules, is given by a hereditary category
$\mathcal{H}$. We refer to~\cite{kussin:07} for more information on
$\mathcal{H}$. This category is determined by \emph{orbit algebras\/}
$$R=\Pi (L,\sigma)=\bigoplus_{n\geq 0}\Hom_{\Lambda} (L,\sigma^n L)$$
where $L$ is a fixed line bundle and $\sigma$ is a positive
automorphism (that is, $\deg (\sigma L)>\deg (L)=0$). If $f\in\Hom
(L,\sigma^n L)$ and $g\in\Hom (L,\sigma^m L)$, then their product
$f\ast g$ in $R$ is defined to be $\sigma^m (f)\circ g$. Note that $R$
is typically noncommutative. With this we have
\begin{equation}
  \label{eq:quotient-cat}
  \mathcal{H}\simeq\mod^{\mathbb{Z}}(R)/\mod_0^{\mathbb{Z}}(R),
\end{equation}
the quotient category formed with respect to the Serre subcategory of
objects of finite length. Hence in the terminology
of~\cite{artinzhang:94}, $\mathcal{H}$ is a (noncommutative)
noetherian projective scheme. The category $\mathcal{H}$ is also
obtained in a more simple minded way by removing the preinjective
component from $\mod (\Lambda)$ and gluing it together with the
preprojective component. The regular modules (lying in homogeneous
tubes) become the objects of finite length in $\mathcal{H}$, with a
simple object $S_x$ for each tube index $x\in\mathbb{X}$. The
Auslander-Reiten translations $\tau$ and $\tau^{-}$ are
autoequivalences of $\mathcal{H}$. The bounded derived categories of
$\mathcal{H}$ and of $\mod\Lambda$ are equivalent, as triangulated
categories. We denote by $L$ an object in $\mathcal{H}$ which
corresponds to a projective module of defect $-1$ and consider it as a
structure sheaf. 

\begin{Defi}[\cite{kussin:07}]
  An automorphism $\sigma $ of $\mathcal{H}$ is called
  \emph{efficient\/}, if it is positive, fixes each tube and $\deg
  (\sigma L)>0$ is minimal with these properties.
\end{Defi}

Note that $\tau^{-}$ and all tubular shifts $\sigma_x$ are positive,
fixing all tubes. But in general minimality is not fulfilled for these
automorphisms. On the other hand efficient automorphisms always exist,
and in many cases there is even a tubular shift automorphism which is
efficient. But there are examples of tame bimodules such that there is
no tubular shift which is efficient (\cite[1.1.13]{kussin:07}; see
also below). The reason why considering efficient automorphisms is the
following.

\begin{Thm}[\cite{kussin:07}]
  Let $\sigma$ be an efficient automorphism. Then the orbit algebra
  $R=\Pi (L,\sigma)$ is a (not necessarily commutative) graded
  factorial domain. Moreover, there is a natural bijection between the
  points of $\mathbb{X}$ and the homogeneous prime ideals of height
  one in $R$, and these are left and right principal ideals, having a
  normal element as a generator.
\end{Thm}

Here we use a graded version of the notion of unique factorization
(=factorial) domains introduced by Chatters and
Jordan~\cite{chattersjordan:86}. We denote by $\pi_x$ a normal
(homogeneous) generator of the prime ideal corresponding to $x$. We
call these elements \emph{prime}. The correspondence works via
universal extensions, see below.

Moreover, the function (skew) field $k(\mathbb{X})$ of $\mathbb{X}$ is
obtained as the quotient division ring (of degree zero fractions) of
$R$. This is an algebraic function skew field in one variable, of
finite dimension over its centre; we denote by $s(\mathbb{X})$ the
square root of this dimension. This function field is of importance in
representation theory since it coincides with the endomorphism ring of
the unique generic $\Lambda$-module.

\begin{Prop}[\cite{kussin:07}]\label{Prop:finite-over-centre}
  Let $\sigma$ be an efficient automorphism. Then the orbit algebra
  $R=\Pi (L,\sigma)$ is finitely generated as module over its
  noetherian centre $C$. Accordingly, the map $P\mapsto P\cap C$ is a
  bijection between the set of homogeneous prime ideals of height one
  in $R$ and the set of those in $C$.
\end{Prop}

For each $x\in\mathbb{X}$ denote by $$f(x)=\frac{1}{\varepsilon}[\Hom
(L,S_x):\End (L)],\ e(x)=[\Hom (L,S_x):\End (S_x)]$$ the
\emph{index\/} and the \emph{multiplicity\/}, respectively, of $x$.
Here, $\varepsilon=1$ if $M$ is a $(2,2)$-bimodule, and
$\varepsilon=2$ if $M$ is a $(1,4)$- or $(4,1)$-bimodule. A point $x$
is called \emph{rational\/}, if $f(x)=1$, and \emph{unirational\/}, if
additionally $e(x)=1$. By~\cite[Prop.~4.2]{lenzingdelapena:99} there
is always a rational point $x$. If $L$ is the structure sheaf and
$x\in\mathbb{X}$, then the universal extension of $L$ with respect to
the tube of index $x$ is given by $0\rightarrow L\rightarrow
L(x)\rightarrow S_x^{e(x)}\rightarrow 0$. This is more generally
defined also for other objects, and is functorial, the assignment
$A\mapsto A(x)$ inducing an autoequivalence of $\mathcal{H}$, denoted
by $\sigma_x$. This is called a tubular mutation, or a tubular shift.
We refer to~\cite{lenzingdelapena:99} (or~\cite{kussin:07}) for more
details on this subject.

If $\sigma$ is efficient, then $L(x)=\sigma^d (L)$ for some positive
integer $d$, and the $S_x$-universal extension above becomes
$$0\rightarrow L\stackrel{\pi_x}\rightarrow \sigma^d (L)\rightarrow
S_x^{e(x)}\rightarrow 0,$$ where $\pi_x$ is a prime element of degree
$d$ in $\Pi (L,\sigma)$.

The advantage of an efficient tubular mutation is given by the
following.

\begin{Prop}
  \textnormal{(1)} Let $\sigma_x$ be a tubular shift associated with a
  tube of index $x$. If $\sigma_x$ is efficient then the prime element
  $\pi_x$ is central in the orbit algebra $R=\Pi (L,\sigma_x)$ and of
  degree one.
  
  \textnormal{(2)} If $\sigma$ is efficient and if there is a central
  element in $\Pi (L,\sigma)$ of degree one, then $\sigma=\sigma_x$ is
  a tubular shift.
\end{Prop}
\begin{proof}
  (1) This follows directly from the definition of the multiplication
  in the orbit algebra and functorial properties of universal
  extensions, see~\cite[1.7.1]{kussin:07}. For (2)
  see~\cite[3.2.13]{kussin:07}. 
\end{proof}

\begin{Rems}
  Let $M$ be a tame bimodule with centre $k$. 
  
  (a) There is always an efficient tubular shift $\sigma_x$ in the
  cases where $k$ is algebraically closed, or real closed, or
  (see~\ref{Prop:finite-field-unirational} below) a finite field.

  (b) There is always an efficient tubular shift $\sigma_x$ if $M$ is
  a non-simple bimodule.

  (c) If $M$ is a simple $(2,2)$-bimodule then there is an efficient
  tubular shift $\sigma_x$ if and only if there is a point $x$ such
  that $f(x)=1$ and $e(x)=2$, or $f(x)=2$ and $e(x)=1$.

  (d) If $M$ is a $(1,4)$- (or $(4,1)$-) bimodule, then there is an
  efficient tubular shift $\sigma_x$ if and only if there is a
  unirational point $x$. A class of examples when this happens is
  given by Lemma~\ref{lem:intermediate-degree-two}.

  (e) Let $k=\mathbb{Q}$ and $\zeta$ be a primitive third root of
  unity. Then the tame bimodule $M={}_{\mathbb{Q}(\sqrt[3]{2})}
  \mathbb{Q}(\sqrt[3]{2},\zeta)_{\mathbb{Q}(\zeta\sqrt[3]{2})}$ does
  not allow an efficient tubular shift (\cite[1.1.13]{kussin:07}).
\end{Rems}

We will often consider the following type of a $(1,4)$-bimodule:
$M={}_k F_F $, where $F$ is a skew field, with $k$ lying in its centre
and of dimension four over $k$.

\begin{Lem}\label{lem:intermediate-degree-two}
  Let $F/k$ be a skew field extension of dimension four, and let $K$ be
  an intermediate field of degree two. Then the tame bimodule $M={}_k
  F_F$ admits a simple regular representation $S_x$ with $\End
  (S_x)\simeq K$, where $x$ is a unirational point.
\end{Lem}
\begin{pf}
  Let $K=k(x)$, with\footnote{Here and in the following we abuse
    notation a little bit: we use the letters $x,\,y,\,\dots$ to
    denote points in $\mathbb{X}$ as well as (base) elements in field
    extensions and also generators in graded algebras; but in such a
    case these different entities are so strongly related to each
    other that it is convenient to use the same symbol.} $x\in
  K\setminus k$ and $x^2 =c_1 x+c_0$, with $c_0$, $c_1 \in k$. Then
  $S_x = (k^2 \otimes F\stackrel{(1,x)}\rightarrow F)$ is a simple
  regular representation with $f(x)=1$. Moreover, the correspondence $
  \begin{pmatrix}
    a & bc_0\\
    b & a+bc_1
  \end{pmatrix}\leftrightarrow a+bx$ gives an isomorphism $\End
  (S_x)\simeq k(x)$. Since $[k(x):k]=2$ we get $e(x)=1$.
\end{pf}

\begin{Prop}\label{Prop:finite-field-unirational}
  Let $k$ be a finite field. There is a unirational point
  $x\in\mathbb{X}$, that is, we have $e(x)=1=f(x)$. Accordingly, the
  corresponding tubular shift $\sigma_x$ is efficient.
\end{Prop}
\begin{pf}
  Any $(2,2)$-bimodule over a finite field is non-simple, and thus
  admits a unirational point (compare~\cite[0.6.2]{kussin:07}). If $M$
  is a $(1,4)$-bimodule over a finite field, hence of the form $M={}_k
  K_K$ where $K/k$ is a field extension of degree four, and there is
  an intermediate field of degree two. Now apply the preceding lemma.
\end{pf}

\section{The ghost group}

We denote by $\Aut (\mathbb{X})$ the group of isomorphism classes of
autoequivalences of $\mathcal{H}$ which fix $L$ (up to isomorphism).
We denote by $\mathcal{G}$ the \emph{ghost group}, which is the
subgroup consisting of (the classes of) those autoequivalences which
leave each tube fixed. It is easy to see that it is equivalent to
require that such an autoequivalence fixes all objects of
$\mathcal{H}$ (up to isomorphism). Representatives of the non-trivial
elements of the ghost group are called \emph{ghosts}. We will show in
this section that the occurrence of ghosts is due to noncommutativity.

Let $R=\Pi (L,\sigma)$, where $\sigma$ is efficient. Denote by $\Aut
(R)$ the group of graded algebra automorphisms of $R$, and by $\Aut_0
(R)$ the subgroup of those automorphisms fixing all homogeneous prime
ideals of $R$. For the general definition of the normal subgroup
$\overline{\Inn}(R)$ we refer to~\cite[3.2.2]{kussin:07}. Here we will
consider only a special case. Let $\alpha\in\Aut (R)$ and
$N\in\mod^{\mathbb{Z}}(R)$. Define a new graded module structure on
$N$ by the formula $x\cdot r\stackrel{def}=x\alpha^{-1}(r)$. As shown
in~\cite[3.2.3]{kussin:07}, using~\eqref{eq:quotient-cat}, this
induces an element $\alpha_{\ast}\in\Aut (\mathbb{X})$, and
$\alpha_{\ast}\simeq 1$ (the identity functor on $\mathcal{H}$) if and
only if $\alpha\in\overline{\Inn}(R)$. Moreover,
$\alpha_{\ast}\in\mathcal{G}$ if and only if $\alpha\in\Aut_0 (R)$. In
the case which we consider here, when additionally $\sigma=\sigma_x$
is a tubular shift, the description of $\overline{\Inn}(R)$ is the
following: this group is then generated by the inner automorphisms
defined in the usual way, and by automorphisms $\varphi_a$ ($a\in R_0$
non-zero, lying in the centre of $R$) of the following form:
$\varphi_a (r)=a^n r$ for homogeneous elements $r$ of degree $n$. If,
for example, $R_0 =k$ then $\overline{\Inn}(R)=k^{\ast}$.


\begin{Thm}\label{Thm:ghost-liftable}
  Let $\sigma_x$ be an efficient tubular shift, let $R=\Pi
  (L,\sigma_x)$. 
  Then for the ghost group we have $\mathcal{G}\simeq\Aut_0
  (R)/\overline{\Inn}(R)$.
\end{Thm}
\begin{proof}
  By the preceding remarks there is an injective homomorphism of
  groups $\Aut_0 (R)/\overline{\Inn}(R)\rightarrow\mathcal{G}$,
  induced by $\alpha \mapsto\alpha_{\ast}$. It is therefore sufficient
  to show that a given $\gamma\in\mathcal{G}$ is liftable to a graded
  algebra automorphism of $R$. Denote by $\mathcal{L}$ the full
  subcategory with objects $\sigma_x^n L$ ($n\in\mathbb{Z}$). Changing
  $\gamma$ by a suitable isomorphism, we can assume that $\gamma
  (X)=X$ for all objects $X\in\mathcal{L}$. Denote by $\pi_x$ the
  (central) prime element corresponding to $x$. Since $\gamma$ fixes
  the tube of index $x$ there is a natural isomorphism
  $\eta\colon\gamma\sigma_x\stackrel{\sim}\rightarrow\sigma_x \gamma$,
  and $\gamma (\sigma_x^n (\pi_x))=\sigma_x^n (\pi_x a_n)$ for some
  $a_n \in R_0^{\ast}$ for all $n\in\mathbb{Z}$
  (see~\cite[0.4.8]{kussin:07}). Since $\sigma_x^n
  (\pi_x)\sigma_x^n(a_n)=\sigma_x^{n+1}(a_n) \sigma_x^n (\pi_x)$ it is
  obvious how to change $\gamma$ again by an isomorphism in such a
  way, that we can assume $\gamma (\sigma_x^n (\pi_x))=\sigma_x^n
  (\pi_x)$ for all $n\in\mathbb{Z}$. We obtain $\eta =1$ on
  $\mathcal{L}$. By restriction, $\gamma$ induces an autoequivalence
  of $\mathcal{L}$, and via the section functor
  $\Gamma=\oplus_{n\in\mathbb{Z}}\Hom (L,\sigma_x^n ?)$ also of the
  full subcategory of $\mod^{\mathbb{Z}}(R)$ consisting of the graded
  modules $R(n)$ ($n\in\mathbb{Z}$). Clearly, $\gamma$ induces a
  bijective graded $k$-linear map $\alpha\colon R\rightarrow R$,
  mapping an element $f\in R_n =\Hom (L,\sigma_x^n L)$ to $\gamma
  (f)$. Since $\gamma\sigma_x =\sigma_x \gamma$ on $\mathcal{L}$, it
  follows easily that $\alpha$ is a graded algebra automorphism. Since
  by construction $\alpha_{\ast}$ and $\gamma$ agree on $\mathcal{L}$,
  we get $\alpha_{\ast}\simeq\gamma$ also on $\mathcal{H}$, and the
  claim follows.
\end{proof}

\begin{Rem}\label{Rem:liftable}
  The proof shows that each $\gamma\in\Aut (\mathbb{X})$ which fixes
  the tube of index $x$ is liftable to a graded automorphism of $R=\Pi
  (L,\sigma_x)$.
\end{Rem}

\begin{Cor}
  Let $\mathbb{X}$ be commutative. Then\/ $\mathcal{G}=1$.
\end{Cor}
\begin{proof}
  If $k(\mathbb{X})$ is commutative we have $e(x)=1$ for all
  $x\in\mathbb{X}$ (by~\cite[4.3.1]{kussin:07}). Let $x$ be a rational
  point. Then the tubular shift $\sigma_x$ is efficient, and $R=\Pi
  (L,\sigma_x)$ is commutative. By~\eqref{eq:quotient-cat}, $R_0$ lies
  in the centre $k$ of $\mathcal{H}$, hence $R_0 =k$ follows. Let
  $\gamma\in\Aut_0 (R)$. If $\pi\in R$ is homogeneous prime, $\gamma
  (\pi)=a\pi$ for some $a\in k^{\ast}$.  Since $R$ is a commutative
  graded factorial domain, $a$ only depends on the degree of $\pi$,
  and $\gamma\in\overline{\Inn}(R)$ follows.
\end{proof}

The converse of the corollary does not hold, as the example
$M=\mathbb{H}\oplus\mathbb{H}$ over the real numbers shows. Here
$\mathcal{G}=1$, but $k(\mathbb{X})=\mathbb{H}(T)$ is not commutative.

\section{The $(2,2)$-case over finite fields}

We study now interesting classes of tame bimodules in more detail.
Before we concentrate on the case of a $(1,4)$-bimodule, we briefly
describe the case of a $(2,2)$-bimodule. Here we restrict to the case
when $k$ is a finite field, since then all $(2,2)$-bimodules can be
described explicitly. For a finite-dimensional field extension $K/k$,
an element $\alpha\in\Gal (K/k)$ and an $(\alpha,1)$-derivation
$\delta$ of $K$ denote by $M=M(K,\alpha,\delta)$ the following (tame)
bimodule: as a left $K$-vector space $M=K\oplus K$. As a right
$K$-vector space we have the rule $(x,y)\cdot a=(xa+y\delta
(a),y\alpha (a))$ (for $x$, $y$, $a\in K$). If $\delta=0$ we also
write $M={}_K K_K \oplus {}_K K_{K^{\alpha}}$.

\begin{Prop}\label{Prop:2-2-finite}
  Let $M$ be a tame bimodule over a finite field of dimension type
  $(2,2)$ with centre $k$.  Then there is a finite field extension
  $K/k$ such that $M={}_K M_K$.  More precisely, if $[K:k]=n$, $\Gal
  (K/k)=\langle\alpha\rangle$, then $M={}_K K_K \oplus {}_K
  K_{K^{\alpha}}$. Accordingly, for a unirational point $x$ the
  corresponding orbit algebra $R=\Pi (L,\sigma_x)$ is isomorphic to
  the skew polynomial algebra $K[X;Y,\alpha]$, graded by total degree.
  Its centre is given by $k[X,Y^n]$, and $s(\mathbb{X})=n$. The
  function field $k(\mathbb{X})$ is isomorphic to $K (T,\alpha)$, and
  its centre is given by $k(T^n)$. The ghost group $\mathcal{G}$ is
  cyclic of order $n$, generated by $\alpha_{\ast}$.
\end{Prop}
\begin{pf}
  It is easy to see that over a finite field each $(2,2)$-bimodule is
  non-simple. By~\cite{ringel:76} $M$ is of the form $M=M(K,\alpha,0)$
  for some automorphism $\alpha$ or $M=M(K,1,\delta)$ for some
  derivation $\delta$. In case of a finite field, obviously
  $\delta=0$. Now, $k$ is the centre of $M(K,\alpha,0)$ if and only if
  $k$ is the fixed field of $\alpha$.  The other assertions follow
  from~\cite[5.3.4]{kussin:07}.
\end{pf}

\section{Some $(1,4)$-cases}

The following theorem applies in particular to any tame bimodule of
dimension type $(1,4)$ over a finite field $k$ with $\ch k\neq 2$, but
there are also many other applications. It
generalizes~\cite[1.7.12]{kussin:07} and improves results
of~\cite{baer:83}.

\begin{Thm}\label{Thm:orbit-algebra-commext}
  Let $k$ be a field. Consider the tower of (commutative) fields
  $$k\subsetneq k(x)\subsetneq k(x,y)=K$$
  such that $x^2 =c_0$ and
  $y^2 =a_0 +a_1 x$ for some $c_0$, $a_0$, $a_1 \in k$. (Then $a_1
  \neq 0$ if and only if $k(x)=k(y^2)$.) Let $M$ be the tame
  bimodule $M={}_k K_K$.\\

  \textnormal{(1)} The simple regular representation $$S_x = (k^2
  \otimes K\stackrel{(1,x)}\longrightarrow K)$$ has endomorphism ring
  $k(x)$
  and defines a unirational point.\\

  \textnormal{(2)} The corresponding orbit algebra $R=\Pi
  (L,\sigma_x)$ is the $k$-algebra on three generators $X$, $Y$ and
  $Z$ with relations
  \begin{gather}
    XY-YX=0,\\
    XZ-ZX=0,\\
    ZY+YZ+a_1 X^2 =0,\\
    Z^2 +c_0 Y^2 -a_0 X^2 =0.
  \end{gather}\\
  If $\ch k =2$ and $a_1 =0$ then $R$ is commutative, otherwise its
  centre is given by $k[X,Y^2]$.\\
  
  \textnormal{(3)} The function field $k(\mathbb{X})$ is isomorphic to
  the quotient division ring of $$k\langle U,V\rangle /(VU+UV+a_1 ,V^2
  +c_0 U^2 -a_0).$$
  If $\ch k=2$ and $a_1 =0$ then $k(\mathbb{X})$ is
  commutative, otherwise its centre is $k(U^2)$ and hence
  $s(\mathbb{X})=2$.\\

  \textnormal{(4)} We have $\Aut (\mathbb{X})\simeq\Gal (K/k)$. Assume
  $\ch k\neq 2$. If $a_1 =0$ then $\Aut (\mathbb{X})$ is the Klein
  four group coinciding with the ghost group, and if $a_1 \neq 0$ then
  $\Aut (\mathbb{X})$ is cyclic. In both cases, the graded algebra
  automorphism induced by sending $X$ to $-X$ (and leaving $Y$ and $Z$
  fixed) induces a ghost of order $2$. If $a_1 \neq 0$ and $a_0 =0$,
  then this is the only ghost.
\end{Thm}
\begin{pf}
  For $n\geq 1$ denote by $P_n$ the representation of rank $1$ (defect
  $-1$) given by $k^{2n-1}\otimes K\stackrel{C_n}\longrightarrow K^n$,
  where $$C_n= \left(\begin{array}{ccccc|cccc}
      1 &  & & & & x &  &  &\\
      & 1 & &  & & y & x &  &\\
      &  & \ddots &  &  & & y & \ddots &\\
     &  &  & 1 & &  &  & \ddots & x\\
     &  &  &  &1 &  &  & & y\\
  \end{array}\right)\in\matring_{n,2n-1}(K).$$
These are the preprojective representations of rank $1$. This can be
shown as in~\cite{baer:83} where the assumption $a_0 \neq 0$ is not
needed, but there is another argument, avoiding also many
calculations: First, it is easy to check that $\End (P_2)=k$. Of
course, all $P_n$ are of rank $1$.  Moreover, it is easy to see that
$\Hom (S_x,P_n)=0$ for all $n\geq 1$. In particular, $P_n$ has no
preinjective summand. There is a short exact sequence
\begin{eqnarray}
  \label{eq:univ-ext}
  0\rightarrow P_n \stackrel{X_n}\rightarrow P_{n+1}
\rightarrow S_x \rightarrow 0
\end{eqnarray} for every $n\geq 1$, where $X_n$ is defined below. It
follows that $x$ is a rational point, and since the simple regular
representation $S_x$ has endomorphism ring $\End (S_x)=k(x)$, we see
that $e(x)=1$.  Thus~\eqref{eq:univ-ext} is the $S_x$-universal
extension of $P_n$. In particular, $P_{n+1}\simeq\sigma_x (P_n)$.
Assume by induction that $P_n$ is preprojective with $\End (P_n)=k$.
Then $\End (P_{n+1})\simeq
\End (P_n)=k$, and thus $P_{n+1}$ is also preprojective of rank $1$.

Note that as objects in $\mathcal{H}$ we have $P_n
\simeq\sigma_x^{n-1}L$ for all $n\geq 1$.

A $k$-basis of $\Hom (P_n,P_{n+1})$ is given by the following three
pairs of matrices $X_n =(\overline{X}_n,\overline{\overline{X}}_n)$,
$Y_n =(\overline{Y}_n,\overline{\overline{Y}}_n)$, $Z_n
=(\overline{Z}_n,\overline{\overline{Z}}_n)$, with $\overline{X}_n$,
$\overline{Y}_n$, $\overline{Z}_n \in\matring_{2n+1,2n-1}(k)$ and
$\overline{\overline{X}}_n$, $\overline{\overline{Y}}_n$,
$\overline{\overline{Z}}_n \in\matring_{n+1,n}(K)$ given by
$$\overline{X}_n={\footnotesize
  \left(\begin{array}{cccc|cccc}
      0 &  &  & 0 &  &  & & \\
      1 &  & &  &  &  &  &\\
      & 1 &   & &  &  &  &\\
      &  & \ddots  &  & &  &  &\\
     &  &  & 1 &  &  &  & \\
    \hline
& & & & 0 & & & 0\\
    & & & & 1 & & & \\
    & & & & & 1 & & \\
    & & & & & & \ddots & \\
& & & & & & & 1 \\
  \end{array}\right)}\, ,\ \ \overline{Y}_n={\footnotesize
  \left(\begin{array}{cccc|cccc}
    1 &  & &  &  &  &  &\\
      & 1 &   & &  &  &  &\\
      &  & \ddots  &  & &  &  &\\
     &  &  & 1 &  &  &  & \\
    0 &  &  & 0 &  &  & & \\
    \hline
    & & & & 1 & & & \\
    & & & & & 1 & & \\
    & & & & & & \ddots & \\
& & & & & & & 1 \\
& & & & 0 & & & 0
  \end{array}\right)}$$ and $\overline{\overline{X}}_n$,
$\overline{\overline{Y}}_n$ are the left upper submatrices,
respectively; moreover 
$$\overline{\overline{Z}}_n ={\footnotesize\left(
    \begin{array}{ccccc}
       \ddots &  & & &\\
       \ddots & -x & & &\\
       & -y & x & &\\
       & -a_1& y & -x & \\
       & & & -y & x\\
       & & & -a_1 & y
    \end{array}\right)}$$ 
and 
$$\overline{Z}_n ={\footnotesize\left(
    \begin{array}{rrrrrr|rrrrrr}
       & \ddots & & & & & c'_0   & & & &   \\
       & & & & & & 0 & -c'_0 & & &\\
       & & -a_1& &  & & a'_0 & 0 & c'_0 & &    \\
       & & & 0 &  & & & -a'_0 & & \ddots  &       \\
       & & & & -a_1 & 0 & & & \ddots & &    \\
       \hline
       & \ddots & & & & &  & & & & &  \\
       & & & & &  & & \ddots & & & &    \\
       & & -1 & & & & & & & & &       \\
       & & & 1 &  & & & & -a_1& &  &      \\
       & & &   & -1  & & & & & 0 &  &        \\
       & & & & & 1  &  & & & & -a_1 & 0   
     \end{array}\right)}$$ with $a'_0 =(-1)^n a_0$ and $c'_0
 =(-1)^{n-1}c_0$. It is easy to check that for each $n\geq 1$ the
 relations
 \begin{gather*}
   X_{n+1}Y_n =Y_{n+1}X_n,\\
    X_{n+1}Z_n =Z_{n+1}Y_n,\\
    Z_{n+1}Y_n = -Y_{n+1}Z_n -a_1 X_{n+1}X_n,\\
    Z_{n+1}Z_n = -c_0 Y_{n+1}Y_n +a_0 X_{n+1}X_n
 \end{gather*}
hold. Moreover, $\Hom (P_n,P_{n+t})$ is generated by all the
``monomials'' in the variables $X$, $Y$, $Z$ of degree $t$.

Since $X_n$ acts centrally, it follows from~\eqref{eq:univ-ext} that
the shift $X_n \mapsto X_{n+1}$, $Y_n \mapsto Y_{n+1}$, $Z_n \mapsto
Z_{n+1}$ coincides on the full subcategory $\mathcal{L}_{+}$ given by
the $P_n$ ($n\geq 1$) with the tubular shift $\sigma_x$ associated
with the tube which contains $S_x$.

Assume $\ch k\neq 2$. The graded algebra automorphism induced by
replacing $X$ by $-X$ is prime fixing, since the centre of $R$ is
given by $k[X,Y^2]$. It is easily checked (on the full subcategory of
$\mathcal{H}$ with objects $L$ and $L(1)$) that the induced functor on
$\mathcal{H}$ is not isomorphic to the identity functor.

Let now $a_1 \neq 0$ and $a_0 =0$. Let $\gamma$ be in $\Aut_0 (R)$. Up
to an element in $\overline{\Inn}(R)$ we can assume that $\gamma
(X)=X$, and then also $\gamma (Y^2)=Y^2$ and $\gamma (Z^2)=Z^2$. A
direct calculation then shows that only $\gamma (Y)=Y$, $\gamma (Z)=Z$
or $\gamma (Y)=-Y$, $\gamma (Z)=-Z$ is possible. This proves the
theorem.
\end{pf}

\begin{Rem}
  An affine version of this graded algebra $R$ (by making $X$
  invertible) can be found in~\cite{dlab:83}.
\end{Rem}

\begin{Rem}
  Denote by $x$, $y$, $z$ the classes of $X$, $Y$, $Z$, respectively.
  The elements $x^2$, $xy$, $y^2$, $zx$, $zy$ form a $k$-basis of
  $R_2$. If $a_1 \neq 0$ then it follows easily that the elements $y$
  and $z$ are not normal, but $y^2$ and $z^2$ are prime\footnote{In
    the ``extreme'' case when $a_0 =0$ and $a_1 =1$, the prime
    elements $y^2$ and $z^2$ even coincide, up to multiplication with
    a unit; in particular, the cokernels of $Y$ and $Z$ are
    isomorphic.}; the only prime (or normal) element of degree one is
  $x$ (up to multiplication with a unit).
\end{Rem}

\begin{Exs}\label{Ex:several-examples}
  We keep the assumptions and notations of the preceding theorem. We
  additionally assume $\ch k\neq 2$ in the following. We consider the
  automorphism group of $\mathbb{X}$ and the ghost group, which is
  always non-trivial in this situation.
  
  (1) $a_1 =0$. In this case $\Aut (\mathbb{X})\simeq\Gal
  (K/k)\simeq\mathbb{V}_4$, the Klein four group, generators induced
  by the graded algebra automorphisms $(X,Y,Z)\mapsto (X,Y,-Z)$ and
  $(X,Y,Z)\mapsto (X,-Y,Z)$, and all these automorphisms are ghost
  automorphisms (see~\cite{kussin:07}).
  
  (2) $a_1 \neq 0$. In this case $k(y^2)=k(x)$. 
  \begin{enumerate}
  \item[(i)] Assume $a_1 =1$ and $a_0 =0$. Then $y=\sqrt[4]{c_0}$,
    $x=y^2 =\sqrt{c_0}$, and
  $$\Pi (L,\sigma_x)=k\spitz{X,Y,Z}/\left(
    \begin{array}{c}
      XY-YX,\ XZ-ZX,\\
      YZ+ZY+X^2,\ Z^2 +c_0 Y^2   
    \end{array}\right).$$
  Denote by $\gamma$ the graded algebra automorphism induced by
  $X\mapsto -X$.
  \begin{enumerate}
  \item[a.] $k=\mathbb{Q}$, $K=k(\sqrt[4]{2})$ non-Galois. In this case
    $\Aut (\mathbb{X})\simeq\mathbb{Z}_2$, which coincides with the
    ghost group, and is generated by $\spitz{\gamma_{\ast}}$.
  \item[b.] $k=\mathbb{Q}(i)$, $K=k(\sqrt[4]{2})$ Galois. Then $\Aut
    (\mathbb{X})\simeq\mathbb{Z}_4 =\spitz{\alpha_{\ast}}$, where
    $\alpha$ is the graded algebra automorphism induced by
    $(X,Y,Z)\mapsto (iX,-Y,Z)$; we have $\gamma=\alpha^2$. Here the
    ghost group is generated again by $\gamma_{\ast}$
    (by~\ref{Thm:orbit-algebra-commext}; it is also easy to see that
    $\alpha_{\ast}$ is not a ghost).
  \end{enumerate}
\item[(ii)] $k$ a finite field, with $p^n$ elements. 
  We have $\Aut (\mathbb{X})\simeq\mathbb{Z}_4$, generated by the
  $n$-th power of the Frobenius automorphism, whose square is the
  ghost $\gamma_{\ast}$ of order two. But this is not a ghost itself.
  For example, let $k=\mathbb{F}_3$ and $K=\mathbb{F}_{3^4}=k(x,y)$
  with $x^2 =2$ and $y^2 =x+1$. Then $\alpha (X)=X$, $\alpha
  (Y)=-Y+Z$, $\alpha (Z)=Y+Z$ defines an element of $\Aut (R)$, whose
  square is $\gamma$ (up to an element of $\overline{\Inn} (R)$). But
  $\alpha$ does not fix the prime ideal generated by $Y^2$. Note also,
  that the rational points correspond to the prime elements given by
  the classes of $X$ (the unirational point) and of $Y^2$, $Y^2 +X^2$,
  $Y^2 +2X^2$, which are products of degree one elements.
  \end{enumerate}
\end{Exs}

\begin{Rem}
  We stress that there is a difference between the calculations
  in~\cite{baer:83} and ours: The change of the roles of the matrices
  $X$ and $Y$ is just notation, but more essential and advantageous is
  our change of the roles of $x$ and $y$ in the matrices $C_n$, as the
  preceding example $K=k(\sqrt[4]{2})$ shows: In~\cite{baer:83} only
  relations are obtained which are not invariant under shift $n\mapsto
  n+1$ and where there is no central variable. Moreover,
  in~\cite{baer:83} generators and relations for the category of
  preprojective representations of defect $-1$ were determined, but
  not the functor induced by $n\mapsto n+1$.
\end{Rem}

\begin{Prop}
  Let the assumptions be as in Theorem~\ref{Thm:orbit-algebra-commext}
  (with $\ch k\neq 2$ or $a_1 \neq 0$) or in
  Proposition~\ref{Prop:2-2-finite}. Let $\gamma\in\Aut (\mathbb{X})$
  such that $\gamma (S_p)\simeq S_p$ for the two points $p=x$ and
  $p=y$ (corresponding to the variables $X$ and $Y$, resp.). Then
  $\gamma\in\mathcal{G}$.
\end{Prop}
\begin{proof}
  By~\ref{Rem:liftable} there is an element $\beta\in\Aut_0 (R)$ such
  that $\beta_{\ast}\simeq\gamma$. In both cases the centre of $R$ is
  given by $k[X,Y^n]$ for some natural number $n$, and $\beta$ fixes
  the prime ideals generated by $X$ and $Y^n$, respectively. Hence,
  after possibly changing $\beta$ up to an element of $\overline{\Inn}
  (R)$ we have that $\beta$ induces the identity on the centre. Now
  the result follows by~\ref{Prop:finite-over-centre}.
\end{proof}

\begin{Thm}\label{Thm:orbit-algebra-skewext}
  Let $k$ be a field. Consider the tower of skew fields
  $$k\subsetneq k(x)\subsetneq k\langle x,\,y\rangle=F$$ such that
  $x^2 =c_0$ and $y^2 =a_0$ for some $c_0$, $a_0 \in k$ and with
  $yx=-xy$. Let
  $M$ be the tame bimodule $M={}_k F_F$.\\

  \textnormal{(1)} The simple regular representation $$S_x = (k^2
  \otimes F\stackrel{(1,x)}\longrightarrow F)$$ has endomorphism ring
  $k(x)$
  and defines a unirational point.\\

  \textnormal{(2)} The corresponding orbit algebra $R=\Pi
  (L,\sigma_x)$ is the $k$-algebra on three generators $X$, $Y$ and
  $Z$ with relations
  \begin{gather}
    XY-YX=0,\\
    XZ-ZX=0,\\
    ZY-YZ=0,\\
    Z^2 -c_0 Y^2 -a_0 X^2 =0.\label{eq:last-equation}
  \end{gather}\\
  In particular, $R$ is commutative.
\end{Thm}
\begin{pf}
  Erase all the minus signs and set $a_1 =0$ in the matrices
  $\overline{Z}_n$ and $\overline{\overline{Z}}_n$ in the 
  proof of Theorem~\ref{Thm:orbit-algebra-commext}. 
  The relations are easily verified.
\end{pf}

\begin{Rem}\label{Rem:char-2}
  There is a version of Theorem~\ref{Thm:orbit-algebra-skewext} which
  also applies to quaternion skew fields $\quatt{a}{b}{k}$ in
  characteristic $2$: Let $k$ be a field of characteristic $2$.
  Consider the tower of skew fields
  $$k\subsetneq k(x)\subsetneq k\langle x,\,y\rangle=F$$ such that
  $x^2 =c_0 +x$ and $y^2 =a_0$ for some $c_0$, $a_0 \in k$ ($a_0 \neq
  0$) and with $xy=y+yx$. Let $M$ be the tame bimodule $M={}_k F_F$.
  In that case the relation~\eqref{eq:last-equation} becomes
  \begin{equation}
    \label{eq:rel-quat-2}
    Z^2 +c_0 Y^2 +a_0 X^2 +YZ=0.
  \end{equation}
  In particular, in this case $R$ is also commutative.
\end{Rem}

\section{The commutative case}

In this section we determine those noncommutative homogeneous curves
$\mathbb{X}$ of genus zero which are actually commutative (that is,
$k(\mathbb{X})$ is commutative). For $\ch k\neq 2$ this was done
in~\cite{kussin:07}. We will prove the following result

\begin{Thm}\label{Thm:char-all-comm}
  Let $M$ be a tame bimodule with centre $k$ and $\mathbb{X}$ the
  corresponding homogeneous exceptional curve. Then $\mathbb{X}$ is a
  commutative curve precisely in the following cases (up to duality):
\begin{itemize}
\item[(i)] $M=k\oplus k$ the Kronecker;
\item[(ii)] $M={}_k F_F$, where $F$ is a skew field of quaternions over
  $k$ (arbitrary characteristic);
\item[(iii)] $M={}_k K_K$, where $K/k$ is a four-dimensional field
  extension not containing a primitive element; $K/k$ is biquadratic,
  i.e.\ there are $a$, $c\in k$ such that $K=k(\sqrt{a},\sqrt{c})$.
\end{itemize}
\end{Thm}

Whereas the first two cases correspond to the Brauer-Severi curves, as
pointed out in~\cite{amitsur:55,crawleyboevey:91}, the last case
yields commutative curves $\mathbb{X}$ of genus zero which are not
Brauer-Severi curves, and happens only in characteristic $2$. Note
that by the theorems in the preceding section $k(\mathbb{X})$ and $\Pi
(L,\sigma_x)$ are determined explicitly in all the cases; in the last
case $k(\mathbb{X})$ is the quotient field of $k[U,V]/(cV^2 +aU^2
+1)$. Note also that this result in particular classifies again the
(commutative) function fields in one variable of genus zero (in the
classical sense). Besides the method of the proof, two aspects are
new: The classification in terms of (tame) bimodules, and the more
general, noncommutative context.

Let $M={}_k F_F$, where $F$ is a four-dimensional skew field extension
of $k$ (with $k\subseteq Z(F)$). Let $\Lambda$ be the corresponding
tame hereditary bimodule algebra. If $\ch k\neq 2$, then
$\Lambda/\rad\Lambda$ is a separable $k$-algebra.

\begin{Thm}\label{Thm:commutative-char}
  Let $k$ be a field and $\mathbb{X}$ be a homogeneous exceptional
  curve with centre $k$ and underlying tame bimodule $M$. The
  following are equivalent:
  \begin{enumerate}
  \item[1.] $k(\mathbb{X})$ is commutative.
  \item[2.] For an (equivalently: for each) efficient automorphism
    $\sigma$ the graded algebra $R=\Pi (L,\sigma)$ is commutative.
  \item[3.] We have $e(x)=1$ for all $x\in\mathbb{X}$.
  \item[4.] We have $e(x)=1$ for all rational $x\in\mathbb{X}$.
  \end{enumerate}
    If $\ch k=2$, assume additionally that $\Lambda/\rad\Lambda$ is a
    separable $k$-algebra. Then these conditions are also equivalent
    to the following:
   \begin{enumerate}
  \item[5.]
    \begin{enumerate}
    \item[a.] $M=k\oplus k$, the Kronecker, or

    \item[b.] $M={}_k F_F$ (or ${}_F\! F_k$) where $F$ is a skew field of
      quaternions over $k$.
    \end{enumerate} 
  \end{enumerate}

  If this holds, then in case~a.\ we have $R=k[X,Y]$ and
  $k(\mathbb{X})=k(T)$, the rational function field in one variable
  over $k$, and in case~b.\
  $$R\simeq\begin{cases}
    k[X,Y,Z]/(-aX^2 -bY^2 +abZ^2) & \ch k\neq 2\ \text{and}\
    F=\quat{a}{b}{k};\\
    k[X,Y,Z]/(aX^2 +bY^2 +YZ+Z^2) & \ch k=2\ \text{and}\
    F=\quatt{a}{b}{k};
  \end{cases}$$ and the function field $k(\mathbb{X})$ is given by the
  quotient field of $k[U,V]/(-aU^2 -bV^2 +ab)$ and $k[U,V]/(U^2
  +UV+bV^2 +a)$, respectively.
\end{Thm}
\begin{proof}
  The equivalence of the first four conditions is proved
  in~\cite[4.3.1,\,4.3.5]{kussin:07}. Assuming these conditions, there
  is an efficient tubular shift automorphism $\sigma_x$, and
  by~\cite[4.3.5]{kussin:07} either $\Pi (L,\sigma_x)\simeq k[X,Y]$,
  or $\Pi (L,\sigma_x)\simeq k[X,Y,Z]/(Q)$, where $Q$ is a
  non-degenerate ternary quadratic form.  In the first case $M$ is the
  Kronecker, in the second it is given by ${}_k F_F$, where $F/k$ is a
  four-dimensional skew field extension. By factoriality, in case $\ch
  k\neq 2$, we have that $Q$ is anisotropic over $k$, and hence
  similar to a form $-aX^2 -bY^2 +abZ^2$. In case $\ch k=2$, by
  tensoring with the separable closure $\overline{k}^s$ of $k$ we get
  by~\cite{crawleyboevey:91} the projective line $\mathbb{P}^1
  (\overline{k}^s)$. It follows that the middle term of the almost
  split sequence starting in $L$ has as endomorphism ring a skew field
  of quaternions over $k$, and~5.\ follows from
  Theorem~\ref{Thm:orbit-algebra-skewext} and Remark~\ref{Rem:char-2}.
\end{proof}

\begin{Exs}\label{exs:pure-insep}
  We exhibit two purely inseparable examples. Let $\mathbb{F}_2$ be
  the field with two elements (or it may be any field of
  characteristic two).
  
  (1) Let $K=\mathbb{F}_2 (u,v)$ be the rational function field in two
  variables over $\mathbb{F}_2$. Let $k=\mathbb{F}_2 (u^2,v^2)$. Then
  $K/k$ is a purely inseparable field extension of degree four. Let
  $M$ be the tame bimodule ${}_k K_K$. Denote by $x$ a unirational
  point.  Then by~\ref{Thm:orbit-algebra-commext}
  $$\Pi (L,\sigma_x)\simeq
  k[X,Y,Z]/(Z^2 +u^2 Y^2 +v^2 X^2).$$ The function field
  $k(\mathbb{X})$ is the quotient field of $k[U,V]/(v^2 U^2 +u^2 V^2
  +1)$. Hence, $\mathbb{X}$ is a commutative curve, but not a
  Brauer-Severi curve. Since $\Gal (K/k)=1$ we have $\Aut
  (\mathbb{X})=1$.
  
  (2) Let $K=\mathbb{F}_2 (u)$ be the rational function field in one
  variable over $\mathbb{F}_2$. Let $k=\mathbb{F}_2 (u^4)$. Again,
  $K/k$ is a purely inseparable field extension of degree four. Here
  the function field is the quotient division ring of $k\langle
  U,V\rangle /(UV+VU+1,\,V^2 +u^4 U^2)$, and hence $\mathbb{X}$ is not
  commutative. Again, $\Aut (\mathbb{X})=1$.
\end{Exs}

As an application of the equivalence of conditions~1.\ and~4.\ in
Theorem~\ref{Thm:commutative-char} we will now see the mathematical
reason why additional commutative curves arise in the inseparable
case: it is because of the absence of a primitive element.

\begin{Prop}
  Let $M={}_k K_K$ where $K/k$ is a four-dimensional field extension.
  Then $\mathbb{X}$ is a commutative curve if and only if there is no
  primitive element for $K/k$. If this is the case, then $K/k$ is a
  biquadratic extension.
\end{Prop}
\begin{proof}
  Any simple regular representation corresponding to a rational point
  $y$ is of the form $S_y =(k^2 \otimes K \stackrel{(1,y)}\rightarrow
  K)$ for some $y\in K\setminus k$. If $y$ is a primitive element, one
  gets $\End (S_y)=k$, hence $e(y)=2$. If $y$ is not primitive, then
  $[k(y):k]=2$ and $e(y)=1$ by~\ref{lem:intermediate-degree-two}. By
  applying Theorem~\ref{Thm:commutative-char} the first part of the
  proposition follows.

  If $K/k$ is inseparable, there are two cases: either $K/k$ is purely
  inseparable, or there is an intermediate field $L$ such that $L/k$
  is separable and $K/L$ is purely inseparable. In the second case
  there is always a primitive element. In the first case, if there is
  no primitive element, the extension is clearly biquadratic.
\end{proof}

\begin{Ex}
  Let $K=\mathbb{F}_4 (u)$ and $k=\mathbb{F}_2 (u^2)$, with $u$
  transcendental. Then $K/k$ is a simple inseparable extension and
  thus the corresponding curve $\mathbb{X}$ is not commutative. This
  example is not covered by Theorem~\ref{Thm:orbit-algebra-commext}.
\end{Ex}

\begin{proof}[Proof of Theorem~\ref{Thm:char-all-comm}]
  Let $M={}_G M_F$ be a tame bimodule with centre $k$ such that
  $\mathbb{X}$ is commutative. By~\ref{Thm:commutative-char} we have
  $e(x)=1$ for all $x\in\mathbb{X}$. In particular, there is a
  unirational point $x$. Hence the corresponding tubular shift
  $\sigma_x$ is efficient. If $M$ is a $(2,2)$-bimodule then $\Pi
  (L,\sigma_x)\simeq k[X,Y]$ by~\cite[4.3.5]{kussin:07}, and hence $M$
  is the Kronecker. Hence assume that $M$ is a $(1,4)$-bimodule. Since
  $R=\Pi (L,\sigma_x)$ is commutative, in particular $R_0=\End (L)$ is
  commutative. Say $G\simeq R_0$, thus $M={}_G F_F$. Moreover, since
  $R=R_0 [R_1]$ and because of~\eqref{eq:quotient-cat} we see that $G$
  lies in the centre of $M$, hence $G=k$. So $M={}_k F_F$. If $F$ is
  not commutative, then $F$ is a skew field of quaternions over $k$. If
  $F$ is commutative then it does not contain a primitive element over
  $k$ by the preceding proposition.
\end{proof}

\begin{Rem}
  Denote by $\mathfrak{m}$ the unique maximal homogeneous left ideal
  of $R=\Pi (L,\sigma_x)$, where $\sigma_x$ is an efficient tubular
  shift. In the commutative case the $\mathfrak{m}$-adic completion
  $\widehat{R}$ is a complete factorial domain of Krull dimension two.
  This follows like in~\cite{lenzing:97} by invoking the completion
  functor $\widehat{\ \cdot\ } \colon\CM^{\mathbb{Z}}(R)\rightarrow
  \CM (\widehat{R})$ studied by Auslander and Reiten. For instance,
  the completion of the graded algebra in
  Example~\ref{exs:pure-insep}~(1) is $k[[X,Y,Z]]/(Z^2 +u^2 Y^2 +v^2
  X^2)$, which is hence factorial.
\end{Rem}

\section{The Auslander-Reiten translation}

Over an arbitrary base field the Auslander-Reiten translation is not
determined combinatorially but depends also on the arithmetics of the
base field, due to the existence of ghosts. In~\cite{kussin:07} it was
shown that for the $(2,2)$-bimodule over the reals $M={}_{\mathbb{C}}
\mathbb{C}_{\mathbb{C}}\oplus
{}_{\mathbb{C}}\mathbb{C}_{\overline{\mathbb{C}}}$ the inverse
Auslander-Reiten translation $\tau^{-}$ is not (a power of) a tubular
shift but a product of two different tubular shifts.  Here we present
a similar example of a $(1,4)$-bimodule.

Let $M$ be the $\mathbb{Q}$-$\mathbb{Q}(\sqrt{2},\sqrt{3})$-bimodule
$\mathbb{Q}(\sqrt{2},\sqrt{3})$. The simple regular representation
$$S_x = (k^2 \otimes K\stackrel{(1,\sqrt{3})}\longrightarrow K)$$ has
endomorphism ring $k(\sqrt{3})$ and defines a unirational point. The
corresponding orbit algebra $R=\Pi (L,\sigma_x)$ is a $k$-algebra on
three generators $X$, $Y$, $Z$ with relations $XY-YX=0$, $XZ-ZX=0$,
$ZY+YZ=0$, $Z^2 +3 Y^2 -2 X^2 =0$.  The homogeneous prime ideals
generated by $X$, $Y$ and $Z$, respectively, define the only
unirational points $x$, $y$ and $z$, respectively. The ghost group
$\mathcal{G}$ is isomorphic to the Klein four group, generated by
$\gamma_y^{\ast}$ and $\gamma_z^{\ast}$, induced by $\gamma_y
(X,Y,Z)=(X,Y,-Z)$ and $\gamma_z (X,Y,Z)=(X,-Y,Z)$, respectively.
Therefore, by Example~\ref{Ex:several-examples}~(1)
(with~\cite[Thm.~3.2.8]{kussin:07}) it follows that the only four
efficient automorphisms (up to isomorphism) are given by
$$\sigma_x,\ \sigma_y =\sigma_x \circ\gamma_y^{\ast},\ \sigma_z
=\sigma_x \circ\gamma_z^{\ast},\ \sigma_y \circ\sigma_z \circ
\sigma_x^{-1} \simeq\sigma_x
\circ\gamma_y^{\ast}\circ\gamma_z^{\ast}.$$

\begin{Prop}
  The functors $\tau^{-}$ and $\sigma_x
  \circ\gamma_y^{\ast}\circ\gamma_z^{\ast}$ are isomorphic.
\end{Prop}
\begin{pf}
  By~\cite{brennerbutler:76} $\tau^-$ coincides with the Coxeter
  functor $C^-$. On objects $C^-$ and $\sigma_x$ act the same way.
  Thus, $\sigma_x \circ C^+$ is a ghost automorphism, that is, one of
  $1$, $\gamma_y^{\ast}$, $\gamma_z^{\ast}$ or
  $\gamma_y^{\ast}\circ\gamma_z^{\ast}$. A direct calculation (by K.\
  Dietrich in his diploma thesis) shows that $\sigma_x \circ C^+$ is
  isomorphic on the full subcategory given by the objects $L$ and
  $L(1)$ to $\gamma_y^{\ast}\circ\gamma_z^{\ast}$, that is, the
  automorphism induced by $X\mapsto -X$.
\end{pf}

All four orbit algebras, formed with respect to these efficient
automorphisms, respectively, are graded factorial. The small
preprojective algebra $\Pi (L,\tau^{-})$ seems to have a clear
disadvantage when compared to the other three: it has no central
element of degree one. In fact, since up to isomorphism the functor
$\gamma_y^{\ast}\circ\gamma_z^{\ast}$ is induced by $(X,Y,Z)\mapsto
(-X,Y,Z)$ we get that $\Pi (L,\tau^{-})$ is generated by $X$, $Y$ and
$Z$ having the relations   
\begin{gather*}
    XY+YX=0,\\
    XZ+ZX=0,\\
    ZY+YZ=0,\\
    Z^2 +3 Y^2 +2 X^2 =0.
  \end{gather*}
  Similar statements hold in the
  examples~\ref{Ex:several-examples}~(2). In all known examples where
  $\tau^{-}$ is computed and where the ghost group is non-trivial,
  $\tau^{-}$ is not a tubular shift.

\begin{Ack}
  This paper was completed during my stay as a Guest Professor at the
  Mathematical Institute of the NTNU in Trondheim. I would like to
  thank the members and other visitors of the Algebra Group there for
  various interesting discussions, and in particular Idun Reiten for
  inviting me. Finally, I am grateful to my partner Gordana Stani\'c
  who was still by my side when my work on this paper began. I
  dedicate this article to her memory.
\end{Ack}

\bibliographystyle{amsplain}

\begin{thebibliography}{1}

\bibitem{amitsur:55} 
S.~A. Amitsur, \emph{Generic splitting fields of central simple
  algebras.\/} Ann.\ of Math.\ 62 (1955), 8--43.

\bibitem{artin:67}
E.~Artin, \emph{Algebraic numbers and algebraic functions.\/} Gordon
and Breach Science Publishers, New York-London-Paris 1967, xiii+349
pp.

\bibitem{artinzhang:94}
M.~Artin and J.~J. Zhang, \emph{Noncommutative projective schemes},
Adv.\ Math.\ 109 (1994), no.~2, 228--287.

\bibitem{baer:83} 
D.~Baer, \emph{Einige homologische Aspekte der
    Darstellungstheorie Artinscher Algebren.\/} Dissertation,
  Universit\"at Paderborn, 1983.

\bibitem{bgl:87} 
D.~Baer, W.\ Geigle and H.\ Lenzing, \emph{The
    preprojective algebra of a tame hereditary artin algebra.\/}
  Comm.\ Algebra 15 (1987), 425--457.

\bibitem{brennerbutler:76} 
S.~Brenner and M.~C.~R.~ Butler, \emph{The
    equivalence of certain functors occurring in the representation
    theory of Artin algebras and species.\/} J.\ London Math.\ Soc.\
  (2) 14 (1976), no.\ 1, 183--187.

\bibitem{buchweitz_eisenbud_herzog:87}
R.-O.\ Buchweitz, D.\ Eisenbud and J.~Herzog, \emph{Cohen-Macaulay
modules on quadrics.\/} Singularities, representation of algebras,
and vector bundles (Lambrecht, 1985), 
Lecture Notes in Math., 1273, Springer, Berlin, 1987, pp.~58--116.

\bibitem{chattersjordan:86}
A.~W.\ Chatters and D.~A.\ Jordan, \emph{Noncommutative unique
  factorisation rings.\/} J.\ London Math.\ Soc.\ (2) 33 (1986),
no.~1, 22--32.

\bibitem{crawleyboevey:91}
W.~W.\ Crawley-Boevey, \emph{Regular modules for tame hereditary
  algebras.\/} Proc.\ London Math.\ Soc.\ (3) 62 (1991), no.~3,
490--508. 

\bibitem{crawleyboevey:91b} 
\bysame, \emph{Tame algebras and generic modules.\/} Proc.\ London
Math.\ Soc.\ (3) 63 (1991), no.~2, 241--265.


\bibitem{dlab:83} 
V.~Dlab, \emph{The regular representations of the
    tame hereditary algebras.\/} S\'eminaire d\`{}Alg\`ebre Paul
  Dubreil et Marie-Paul Malliavin, Proceedings, Paris 1982 (35\`eme
  Ann\'ee) (Berlin-Heidelberg-New York) (M.-P.\ Malliavin, ed.),
  Lecture Notes in Math., vol.\ 1029, Springer-Verlag, 1983,
  pp.~120--133. 

\bibitem{dlabringel:76}
V.\ Dlab and C.\ M.\ Ringel, \emph{Indecomposable representations of
  graphs and algebras.\/} Mem.\ Amer.\ Math.\ Soc.~6 (1976), no.~173,
v+57.   

\bibitem{geiglelenzing:91}
W.~Geigle and H.~Lenzing, \emph{Perpendicular categories with
  applications to representations and sheaves.\/} J.\ Algebra~144
(1991), 273--343. 

\bibitem{kussin:07} D.~Kussin, \emph{Noncommutative curves of genus
    zero -- related to finite dimensional algebras.\/} Mem.\ Amer.\ 
  Math.\ Soc.\ (in press)

\bibitem{lang_tate:52}
S.\ Lang and J.~Tate,
\emph{On Chevalley's proof of Luroth's theorem.\/}
Proc. Amer. Math. Soc. 3, (1952), 621--624.

\bibitem{lenzing:97} 
H.~Lenzing, \emph{Representations of finite
    dimensional algebras and singularity theory.\/} Trends in ring
  theory (Miskolc, 1996) (V.~Dlab et~al., ed.), CMS Conf.\ Proc.,
  vol.~22, Amer.\ Math.\ Soc., Providence, RI, 1998, pp.~71--97.

\bibitem{lenzingdelapena:99} 
H.~Lenzing and J.~A. de~la Pe{\~{n}}a,
  \emph{Concealed-canonical algebras and separating tubular
    families\/}, Proc.\ London Math.\ Soc.\ (3) 78 (1999), no.~3,
  513--540.

\bibitem{ringel:76}
C.~M.~Ringel, \emph{Representations of $K$-species and bimodules.\/}
J.~Algebra 41 (1976), no.\ 2, 269--302.

\bibitem{ringel:79}
\bysame, \emph{Infinite dimensional representations of finite
  dimensional hereditary algebras.\/} Symposia Mathematica, Vol.\
XXIII (Conf.\ Abelian Groups and their Relationship to the Theory of
Modules, INDAM, Rome, 1977), vol.~23, Academic Press, London, 1979,
pp.\ 321--412.

\end{thebibliography}

\providecommand{\bysame}{\leavevmode\hbox to3em{\hrulefill}\thinspace}
\providecommand{\MR}{\relax\ifhmode\unskip\space\fi MR }
\providecommand{\MRhref}[2]{%
  \href{http://www.ams.org/mathscinet-getitem?mr=#1}{#2}
}
\providecommand{\href}[2]{#2}

\end{document}